\documentclass[12pt]{amsart}
\usepackage{amscd,amsthm,amsfonts,amssymb,amsmath,euscript}
\usepackage[dvips]{graphicx}

\voffset= - 0.5in
\hoffset= - 0.3in         

\begin{document}
\begin{flushright}
\hfill FIAN/TD-25/09\\
\hfill ITEP/TH-37/09
\end{flushright}

\bigskip

\centerline{\Large {\bf Universal Algebras of Hurwitz Numbers}}

{\ }

\centerline {\bf A.Mironov\footnote{P.N.Lebedev Physical Institute and ITEP; mironov@itep.ru},
A.Morozov\footnote{ITEP; morozov@itep.ru} and S.Natanzon\footnote{Moscow State University, ITEP
and Moscow Independent University; natanzons@mail.ru}}

{\ }

\centerline{ABSTRACT}

\bigskip

{\footnotesize Infinite-dimensional universal
Cardy-Frobenius algebra is constructed,
which unifies all particular algebras of closed
and open Hurwitz numbers and
is closely related to the algebra of differential operators,
familiar from the theory
of Generalized Kontsevich Model.}

\bigskip


{\bf 1. Introduction.} Classical Hurwitz numbers of complex
algebraic curves generate a commutative Frobenius algebra $A_m$,
which is naturally isomorphic to the center of the group algebra of
symmetric group $\mathbb{S}_m$ \cite{D,MMN}. A natural extension is
provided by Hurwitz numbers of seamed surfaces or foams
\cite{AN2,AN3}. These numbers determine a non-commutative Frobenius
algebra $B_m$ and a homomorphism $\phi_m: A_m\rightarrow B_m$. This
set of data forms a Cardy-Frobenius algebra, which describes Klein
topological field theories \cite{AN1,AN4}. In the present paper the
infinite-dimensional algebras  $A$, $B$ are described which unify
all the Hurwitz numbers algebras.
The construction is based on representation of the group $\mathbb{S}_\infty$
in the algebra ${\mathcal M}$ of formal differential operators, made from the
matrix elements or directly from  $gl(\infty)$ generators,
which has its own value.
${\mathcal M}$ is actually the regular representation of the universal enveloping
algebra of $gl(\infty)$ and
our $B$ is a subset in ${\mathcal M}$,
obtained by taking a kind of operator "traces".
The homomorphism $\phi: A\rightarrow {\mathcal M}$
coincides with the representation of "cut-and-join" operators
constructed in \cite{MMN},
which naturally appear in the theory of Kontsevich integrals
\cite{Ko}-\cite{MSh} and form an associative algebra, isomorphic to the algebra
of ${\mathbb S}_\infty$ characters introduced in \cite{MMN}.
The simplest operators from $B$ (in a different form) appeared in \cite{N}.
It would be also interesting to find a place for the algebras from
\cite{LN} in this context.

\bigskip

{\bf 2. Operators.}
Let $\mathbb{N}$ denote the set of natural numbers (positive integers).
Let $D_{ab}$ with $a,b\in \mathbb{N}$ be $gl(\infty)$ generators. They satisfy
the commutation relation $[D_{ab},D_{cd}] = \delta_{bc}D_{ad} - \delta_{ad}D_{cb}$
and can be conveniently represented in the regular representation
by the differential operators
$D_{ab} = \sum_{e=1}^N X_{ae}\frac{\partial}{\partial X_{be}}$,
and we denote $\tilde{\mathbb{N}}=\{1\ldots N\}\subset\mathbb{N}$.
Introduce "balanced" operators
$$\mathcal{V}_{a_1\ldots\,a_m|\,b_1\ldots\,b_m}
= \ :\!D_{a_1b_1}\ldots D_{a_mb_m}\!:\ =
\!\!
\sum_{(e_1 \ldots\, e_m)\in\, \tilde{\mathbb{N}}^m} X_{a_1e_1}\ldots X_{a_me_m}
\frac{\partial}{\partial X_{b_1e_1}}\ldots \frac{\partial}{\partial X_{b_me_m}}$$
which form a basis in the universal enveloping algebra $Ugl(\infty)$.
The second part of the formula is the explicit definition of the normal ordering in the
first part.
"Balanced" means that the number of $X$'s is the same as
the number of "momenta" $\partial/\partial X$.
The algebra ${\mathcal M}$ consists of linear combinations of such balanced operators.
The algebras $A$ and $B$ are formed by summation over free indices $a_1,\ldots,b_m$
in two different ways.

\bigskip

{\bf 3. Algebra $A$.}
The group $\mathbb{S}_\infty$ is formed by the permutations in $\mathbb{N}$,
involving only the finite sets of numbers.
Define a representation $\phi: \mathbb{S}_\infty \rightarrow {\mathcal M}$,
mapping $\sigma \in \mathbb{S}_m$ into the sum
$$W[\Delta_\sigma]=\sum_{(a_1\ldots\,a_m)\in\tilde{\mathbb{N}}^m}
\mathcal{V}_{a_1\ldots\,a_m|a_{\sigma (1)}\ldots\,a_{\sigma (m)}}$$
This operator depends only on the conjugation class of $\sigma$
(Young diagram $\Delta_\sigma$)
and it is the cut-and-join operator constructed in \cite{MMN}.
These operators form an associative commutative algebra $A \subset {\mathcal M}$.

\bigskip

{\bf 4. Two-fold graphs and algebra $B$.}
A graph $(V_a,E,V_B)$ is called {\it two-fold}, if its vertices are divided
into two sets $V_a$ and $V_b$, and all edges from $E$ have one end in $V_a$
and another in $V_b$.
A homeomorphism of graphs
$\varphi:(V_a,E,V_b)\rightarrow(V'_a,E',V'_b)$
is called isomorphism if   $\varphi(V_a)=(V'_a)$ and $\varphi(V_b)=(V'_b)$.
Denote through $[(V_a,E,V_b)]$ the isomorphism class of $(V_a,E,V_b)$.

Let $\Gamma=[(V_a,E,V_b)]$.
Associate with every edge $E_i \in E$ a pair of numbers $(a_i,b_i)$ so that
$a_i=a_j$, iff  $E_i$ and $E_j$ have a common vertex in $V_a$,
while $b_i=b_j$ iff $E_i$ and $E_j$ have a common vertex in $V_b$.
We call the corresponding operator $\mathcal{V}_{a_1...\,a_m|b_1...\,b_m}$
with $m=\#(E)$ {\it compatible} with  $\Gamma$.
Then, as a straightforward generalization of the above definition of $W[\sigma]$,
denote through $\mathcal{V}(\Gamma)\subset\mathcal{M}$ a sum over
$\{ a_1,\ldots,b_m\}$ of all operators,
compatible with $\Gamma$, with certain normalization factor $c(\Gamma|N)$.
Denote through $B$ the associative but non-commutative algebra formed
by all the operators $\mathcal{V}(\Gamma)$.

Let $\mathcal{B}_m$ be the set of isomorphism classes of the two-fold graphs
with $m$ edges. The associated vector space $B_m$ has a natural structure of
Frobenius algebra, see s.2.3 of \cite{AN3}.
The product of classes  $\Gamma^1=[(V^1_a,E^1,V^1_b)]$ and
$\Gamma^2=[(V^2_a,E^2,V^2_b)]$ is a linear combination of classes $\Gamma$,
consisting of graphs of the form $(V^1_a,E,V^2_b)$
obtained by identification of vertices of the same valence from $V^1_b$ and $V^2_a$
and gluing together the attached edges from $E^1$ and $E^2$.
The structure constants in
$\Gamma^1\Gamma^2=\sum_{\Gamma} C^{\Gamma}_{\Gamma^1\Gamma^2}\Gamma$
take graphs automorphisms into account.
As generalization of a similar statement \cite{MMN} for $A$ we have:

\bigskip

\noindent {\bf Theorem.} \textit{
The structure constants of the algebra $B$,
$\ \mathcal{V}(\Gamma^1)\mathcal{V}(\Gamma^2)=\sum_{\Gamma\in\mathcal{B}}
C^{\mathcal{V}(\Gamma)}_{\mathcal{V}(\Gamma^1)\mathcal{V}(\Gamma^2)}
\mathcal{V}(\Gamma),$
contain the structure constants of all $B_m$
in the following sense:
$\lim_{N\rightarrow\infty}
C^{\mathcal{V}(\Gamma)}_{\mathcal{V}(\Gamma^1)\mathcal{V}(\Gamma^2)}=
C^{\Gamma}_{\Gamma^1\Gamma^2}$, provided $|\Gamma^1|=|\Gamma^2|=|\Gamma|=m$.
The structure constants of $A$ are independent of $N$.
}
\bigskip

{\bf 5. Acknowledgment.} Our work is partly supported by Russian Federal Nuclear Energy
Agency, by RFBR grants 07-02-00878 (A.Mir.), 07-02-00645 (A.Mor.), 07-01-00593 (S.N.),
by joint grants 09-02-90493-Ukr,
09-02-93105-CNRSL, 09-01-92440-CE, 09-02-91005-ANF and by Russian President's Grants
of Support for the Scientific Schools NSh-3035.2008.2 (A.M.'s) and NSh-709.2008.1 (S.N.)


\begin{thebibliography}{12}

\bibitem{D} R.Dijkgraaf,
\emph{Mirror symmetry and elliptic curves, The moduli spaces of curves},
Prog.in Math., {\bf 129} (1995) 149-163

\bibitem{MMN} A.Mironov, A.Morozov and S.Natanzon,
\emph{Complete set of cut-and-join operators in Hurwitz-Kontse\-vich theory},
arXiv:0904.4227

\bibitem{AN2} A.Alexeevski and S.Natanzon,
\emph{ Algebra of Hurwitz numbers for seamed surfaces},
Russian Math.Surveys  61 [4] (2006) 767-769

\bibitem{AN3} A.Alexeevski and S.Natanzon,
\emph{Algebra of two-fold graphs and Hurwitz numbers for seamed surfaces},
72 [4] (2008) 3-24

\bibitem{AN1} A.Alexeevski and S.Natanzon,
\emph{Noncommutative two-dimensional
topological field theories and Hurwitz numbers for real algebraic curves},
Selecta Math., New ser., 12 [3] (2006) 307-377, math.GT/0202164

\bibitem{AN4} A.Alexeevski and S.Natanzon,
\emph{Hurwitz numbers for regular coverings of surfaces by seamed surfaces
and Cardy-Frobenius algebras of finite groups},
Amer.Math.Soc.Transl. 224 [2] (2008) 1-25

\bibitem{Ko} M.Kontsevich,
\emph{Intersection theory on the moduli space of curves and the Airy function},
Comm.Math.Phys. {\bf 147} (1992) 1-23

\bibitem{GKM} S.Kharchev, A.Marshakov, A.Mironov, A.Morozov and A.Zabrodin,
\emph{Towards unified theory of $2d$ gravity},
Nucl. Phys. \textbf{B380} (1992) 181-240, hep-th/9201013

\bibitem{UFN3} A.Morozov,
\emph{Integrability and Matrix Models},
Phys.Usp. \textbf{37} (1994) 1-55, hep-th/9303139;
\emph{Matrix Models as Integrable Systems},
hep-th/9502091

\bibitem{Mirr} A.Mironov,
2d gravity and matrix models. I. 2d gravity,
Int.J.Mod.Phys. {\bf A9} (1994) 4355, hep-th/9312212

\bibitem{AMMP} A.Alexandrov, A.Mironov, A.Morozov and P.Putrov,
\emph{Partition Functions of Matrix Models as the First Special Functions of String Theory.
II. Kontsevich Model},
arXiv:0811.2825

\bibitem{MSh} A.Morozov and Sh.Shakirov,
\emph{Generation of Matrix Models by W-operators},
JHEP 0904:064, 2009, arXiv:0902.2627

\bibitem{N}  S.Natanzon,
\emph{Disk single Hurwitz numbers},
to appear in Funk. An. and Apl., arXiv:0804.0242)

\bibitem{LN} S.Loktev and S.Natanzon,
\emph{Generalized Topological Field Theories from Group Representations},
arXiv:0910.3813


\end{thebibliography}
\end{document}